# POISSOIN ARRIVALS, EXPONENTIAL SERVICE TIME, AND INFINITE SERVERS QUEUE BUSY PERIOD AND BUSY CYCLE DISTRIBUTION FUNCTIONS BOUNDS

**MANUEL ALBERTO M. FERREIRA**
Instituto Universitário de Lisboa (ISCTE – IUL), ISTAR - IUL, Lisboa, Portugal
manuel.ferreira@iscte.pt

**ABSTRACT**

The busy period length distribution function knowledge is important for any queue system, and for the $M|G|\infty$ queue. But the mathematical expressions are in general very complicated, with a few exceptions, involving usually infinite sums and multiple convolutions. So, in this work are deduced some bounds for the $M|M|\infty$ system busy period length distribution function, meaning the second $M$ exponential service time, which analytic expressions are simpler than the exact one. As a consequence, also some bounds for the $M|M|\infty$ system busy cycle length distribution function are presented.

## 1. INTRODUCTION

In a $M|G|\infty$ queue system, $\lambda$ is the Poisson process arrivals rate, $\alpha$ is the mean service time, $G(.)$ is the service time distribution function and so $\alpha = \int_0^\infty [1 - G(t)] \, dt$. $F(.)$ is the service time equilibrium distribution function which expression is $F(t) = \frac{1}{\alpha} \int_0^t [1 - G(x)] dx$. The traffic intensity is $\rho = \lambda\alpha$ and $B$ is the busy period length.

Note the importance of the busy period study, for this queuing system, because in it any customer, when arrives, must find immediately an available server. So, the problem is "for how long the servers – and how many servers – must be available? That is: how long is the busy period length?"

The $B$ distribution function has not a simple form and it can be written as, see (1),

$$P(B \leq t) = 1 - \lambda^{-1} \sum_{n=1}^\infty c^{*n}(t), \ t \geq 0 \qquad (1.1)$$

where $c^{*n}$ is the $n^{th}$ convolution of c with itself being

$$c(t) = \lambda(1 - G(t))e^{-\lambda \int_0^t [1-G(x)]dx} \qquad (1.2).$$

Only for the service time distribution function given by[1], see (2, 3),

$$G(t) = 1 - \frac{(1-e^{-\rho})(\lambda+\beta)}{\lambda e^{-\rho}(e^{(\lambda+\beta)t}-1)+\lambda}, t \geq 0, -\lambda \leq \beta \leq \frac{\lambda}{e^\rho - 1} \quad (1.3)$$

the expression (1.1) becomes simple[2]:

$$P(B \leq t) = 1 - \frac{(1-e^{-\rho})(\lambda+\beta)}{\lambda} e^{-\rho(\lambda+\beta)t}, t \geq 0, -\lambda \leq \beta \leq \frac{\lambda}{e^\rho-1} \quad (1.4).$$

This does not happen for the $M|M|\infty$ queuing systems – exponential service times – and so, in this work, some simple bounds will be presented for $P(B^M \leq t)$.

An idle period followed by a busy period is a busy cycle. Calling Z the busy cycle length,

$$Z = I + B \quad (1.5),$$

where I is the idle period length. So, bounds for $P(Z^M \leq t)$, obtained with simple manipulations after the $P(B^M \leq t)$ ones, and much simpler than the exact expression, will also be presented.

---

[1] It results from, see (4),

$$G(t) = 1 - \frac{1}{\lambda} \frac{(1-e^{-\rho})e^{-\lambda t - \int_0^t \beta(u)du}}{\int_0^\infty e^{-\lambda w - \int_0^w \beta(u)du} dw - (1-e^{-\rho})\int_0^t e^{-\lambda w - \int_0^w \beta(u)du} dw}, t \geq 0, -\lambda \leq \frac{\int_0^t \beta(u)du}{t} \leq \frac{\lambda}{e^\rho - 1},$$

making $\beta(t) = \beta$ (constant).

[2] It results from, see still (4),

$$P(B \leq t) = \left(1 - (1-G(0))\left(e^{-\lambda t - \int_0^t \beta(u)du} + \lambda \int_0^t e^{-\lambda w - \int_0^w \beta(u)du} dw\right)\right) * \sum_{n=0}^\infty \lambda^n (1-G(0))^n \left(e^{-\lambda t - \int_0^t \beta(u)du}\right)^{*n},$$

$$-\lambda \leq \frac{\int_0^t \beta(u)du}{t} \leq \frac{\lambda}{e^\rho - 1}$$

making $\beta(t) = \beta$ (constant).

Finally note that if the probability distribution function allows the study of the distribution structure, only the distribution function allows the probabilities direct calculation.

## 2. BOUNDS FOR $P(B^M \leq t)$

Write c (t) as

$$c(t) = \rho f(t) e^{-\rho F(t)} \qquad (2.1)$$

where $f(t) = \frac{dF(t)}{dt}$. So,

$$c(t) \geq \rho f(t) e^{-\rho} \qquad (2.2)$$

and

$$P(B \leq t) \leq 1 - \lambda^{-1} \sum_{n=1}^{\infty} f^{*n}(t) \rho^n e^{-n\rho}, t \geq 0 \qquad (2.3)$$

or

$$c(t) \leq \rho f(t) \qquad (2.4)$$

and

$$P(B \leq t) \geq 1 - \lambda^{-1} \sum_{n=1}^{\infty} f^{*n}(t) \rho^n, t \geq 0 \qquad (2.5).$$

The $n^{th}$ convolution of $f$ with itself, $f^{*n}$, is the probability density function of the sum of $n$ independent and identically distributed random variables which distribution function is given by $F(t)$. Then the bounds given in (2.3) and (2.5) depend only on $\rho$, $\lambda$ and $F(.)$.

For the $M|M|\infty$ queue, $G(t) = 1 - e^{-\frac{t}{\alpha}}$ and, so, $f(t) = \frac{1}{\alpha} e^{-\frac{t}{\alpha}}$. Then,

$f^{*n}(t) = \left(\frac{1}{\alpha} e^{-\frac{t}{\alpha}}\right)^{*n} = \frac{t^{n-1}}{\alpha^n (n-1)!} e^{-\frac{t}{\alpha}}$. And $\sum_{n=1}^{\infty} f^{*n}(t) \rho^n e^{-n\rho} =$

$\sum_{n=1}^{\infty} \frac{t^{n-1}}{\alpha^n (n-1)!} e^{-\frac{t}{\alpha}} \rho^n e^{-n\rho} = \frac{\rho}{\alpha} e^{-\rho} e^{-\frac{t}{\alpha}} \sum_{n=1}^{\infty} \frac{1}{(n-1)!} \left(\frac{\rho}{\alpha} e^{-\rho} t\right)^{n-1} = \lambda e^{-\rho - \frac{t}{\alpha}} e^{\lambda e^{-\rho} t} =$

$\lambda e^{-\rho + \left(\lambda e^{-\rho} + \frac{1}{\alpha}\right) t} = \lambda e^{-\rho + \frac{\rho e^{-\rho} - 1}{\alpha} t}$.

So

$$P(B^M \leq t) \leq 1 - e^{-\rho - \frac{1-\rho e^{-\rho}}{\alpha} t}, \qquad t \geq 0 \qquad (2.6)$$

after (2.3).

From (2.5), as $\sum_{n=1}^{\infty} f^{*n}(t)\rho^n = \sum_{n=1}^{\infty} \frac{t^{n-1}}{\alpha^n(n-1)!} e^{-\frac{t}{\alpha}} \rho^n =$
$\frac{\rho}{\alpha} e^{-\frac{t}{\alpha}} \sum_{n=1}^{\infty} \frac{1}{(n-1)!} \left(\frac{\rho}{\alpha} t\right)^{n-1} = \lambda e^{-\frac{t}{\alpha}} e^{-\lambda t}$ it is concluded that

$$P(B^M \leq t) \geq 1 - e^{-\frac{1-\rho}{\alpha}t}, \quad t \geq 0 \qquad (2.7).$$

The bound given by (2.6) is always lesser than 1. The one given by (2.7) is positive only for $\rho < 1$.

Otherwise $1 - e^{-\rho - \frac{1-\rho e^{-\rho}}{\alpha}t} \geq 1 - e^{-\frac{1-\rho}{\alpha}t} \Leftrightarrow -\rho - \frac{1-\rho e^{-\rho}}{\alpha}t \leq -\frac{1-\rho}{\alpha}t \Leftrightarrow$
$\frac{1-\rho-1+\rho e^{-\rho}}{\alpha}t \leq \rho \Leftrightarrow \frac{\rho(e^{-\rho}-1)}{\alpha}t \leq \rho \Leftrightarrow t \geq \frac{\rho\alpha}{\rho(e^{-\rho}-1)} = -\frac{\alpha}{(1-e^{-\rho})} < 0, \rho > 0$ and the bound given by (2.6) is always greater than the one given by (2.7).

In (5) it was proved that

$$G(t)e^{-\rho} \leq P(B \leq t) \leq G(t), t \geq 0 \qquad (2.8).$$

Consequently

$$\left(1 - e^{-\frac{t}{\alpha}}\right)e^{-\rho} \leq P(B^M \leq t) \leq 1 - e^{-\frac{t}{\alpha}}, \quad t \geq 0 \qquad (2.9)$$

The lower bound given in (2.9) is always positive, but for $\rho < 1$ the one given by (2.7) is better.

As $1 - e^{-\rho - \frac{1-\rho e^{-\rho}}{\alpha}t} \leq 1 - e^{-\frac{t}{\alpha}} \Leftrightarrow -\rho - \frac{1-\rho e^{-\rho}}{\alpha}t \geq -\frac{t}{\alpha} \Leftrightarrow \frac{1-1+\rho e^{-\rho}}{\alpha}t \geq \rho \Leftrightarrow$
$t \geq \alpha e^\rho$ the upper bound given by (2.6) is better than the one given by (2.9) if $t \geq \alpha e^\rho$.

### 3. BOUNDS FOR $P(Z^M \leq t)$

The random variable $I$ is exponentially distributed with parameter $\lambda$, as it is the case of any queue system with Poisson arrivals. In (6) it was shown that I and B are independent. So, the distribution of Z is the convolution of those of $I$ and $B$.

In fact, the Z distribution function has not a simple form except for the service time given by (1.3), see (7)[3],

---

[3] -For $\beta = 0, P(Z \leq t) = 1 - e^{-e^{-\rho}\lambda t}, t \geq 0$. So Z is exponentially distributed at rate $e^{-\rho}\lambda$ and it may be concluded that the points in time at which begin busy cycles, occur according to a Poisson Process at rate $e^{-\rho}\lambda$.

-For $\beta = \frac{\lambda}{e^\rho-1}, P(Z \leq t) = 1 - \frac{e^\rho-1}{e^\rho-2} e^{-\frac{\lambda}{e^\rho-1}t} + \frac{1}{e^\rho-2} e^{-\lambda t}, t \geq 0$, if $\rho \neq \ln 2$. For $\rho = \ln 2$, using l'Hospital's rule it is obtained $P(Z \leq t) = 1 - (1 + \lambda t)e^{-\lambda t}, t \geq 0$.

$$P(Z \le t) = 1 - \frac{(1-e^{-\rho})(\lambda+\beta)}{\lambda - e^{-\rho}(\lambda+\beta)} e^{-e^{-\rho}(\lambda+\beta)t} + \frac{\beta}{\lambda - e^{-\rho}(\lambda+\beta)} e^{-\lambda t},$$
$$t \ge 0, -\lambda \le \beta \le \frac{\lambda}{e^\rho - 1} \quad (3.1).$$

This does not happen for the $M|M|\infty$ queuing systems and so there will be givens some simple bounds for $P(Z^M \le t)$.

After (2.6) and (2.7), and performing the adequate convolutions

$$1 - \frac{(\rho-1)e^{-\lambda t} + \rho e^{-\frac{1-\rho}{\alpha}t}}{2\rho - 1} \le P(Z^M \le t) \le 1 - \frac{(\rho-1)e^{-\lambda t} + \rho e^{-\rho \cdot \frac{1-\rho e^{-\rho}}{\alpha}t}}{\rho(1+e^{-\rho}) - 1},$$

$$t \ge 0 \quad (3.2).$$

After (2.9), with the adequate convolutions,

$$e^{-\rho}\left(1 - \frac{\rho e^{-\frac{t}{\alpha}} - e^{-\lambda t}}{\rho - 1}\right) \le P(Z^M \le t) \le 1 - \frac{\rho e^{-\frac{t}{\alpha}} - e^{-\lambda t}}{\rho - 1}, t \ge 0 \quad (3.3).$$

## 4. CONCLUSIONS

Upper and lower bounds for $P(B \le t)$ that can be used for every service time distribution were presented. But, for exponential service times, they originate very simple expressions. And it is even possible to compare them to make the best option in their use through very simple rules.

Although the busy cycle is not so important as the busy period its study is of great interest. Then this study was finished with the presentation of bounds for $P(Z^M \le t)$. They are also a good alternative to the exact expression that is not practical at all.